\documentclass{tran-l}
\usepackage{amsmath,amssymb,amscd,amsthm}
\usepackage[matrix,arrow]{xy}

\newcommand{\Z}         {\mathbb Z}

\newcommand{\R}         {\mathbb R}

\newcommand{\BB}        {\mathcal B}
\newcommand{\GG}        {\mathcal G}
\newcommand{\GGo}       {\GG_{0}}
\newcommand{\GGl}       {\GG_{1}}
\newcommand{\HHH}       {\mathcal H}
\newcommand{\HHHo}       {\HHH_{0}}
\newcommand{\HHHl}       {\HHH_{1}}
\newcommand{\KK}        {\mathcal K}
\newcommand{\MM}        {\mathcal M}

\newcommand{\UU}        {\mathcal U}

\newcommand{\Desc}      {\mathrm{Desc}}

\newcommand{\la}        {\leftarrow}
\newcommand{\ra}        {\rightarrow}
\newcommand{\lra}       {\longrightarrow}

\newcommand{\intersect}{\cap}
\newcommand{\iso}       {\cong}
\newcommand{\Prod}      {\displaystyle \prod}

\newcommand{\Geff}      {{\mathcal G}_{\mathrm{eff}}}

\newcommand{\Gtop}      {{\mathcal G}_{\mathrm{top}}}
\newcommand{\Htop}      {{\mathcal H}_{\mathrm{top}}}

\DeclareMathOperator{\Ob}{Ob}
\DeclareMathOperator{\Ar}{Ar}
\DeclareMathOperator{\Hom}{Hom}
\DeclareMathOperator{\Aut}{Aut}
\DeclareMathOperator{\HS}{HS}

  \newtheorem{theorem}{Theorem}[section]
  \newtheorem{Def}[theorem]{Definition}
  \newtheorem{proposition}[theorem]{Proposition}
  \newtheorem{cor}[theorem]{Corollary}
  \newtheorem{lemma}[theorem]{Lemma}

\theoremstyle{definition}
  \newtheorem{DefProp}[theorem]{Definition-Proposition}

\numberwithin{equation}{section}

\begin{document}

\title{Presentations of Noneffective Orbifolds}

\author{Andre Henriques}
\address{Department of Mathematics, University of California, Berkeley CA}
\email{andrhenr@math.berkeley.edu}

\author{David S. Metzler}
\address{Department of Mathematics, University of Florida, Gainesville FL}
\email{metzler@math.ufl.edu}


\date{February 12, 2003}


\begin{abstract}
  It is well-known that an effective orbifold $M$ (one for which the
  local stabilizer groups act effectively) can be presented
  as a quotient of a smooth manifold $P$ by a locally free action
  of a compact lie group $K$. We use the language of groupoids
  to provide a partial answer to the question of whether a
  noneffective orbifold can be so presented. We also
  note some connections to stacks and gerbes.
\end{abstract}

\maketitle

\section{Introduction}
\label{sec:intro}

Recently, motivated largely by string theory, many researchers
have developed new aspects of the geometry and topology of
orbifolds (e.g. \cite{Ruan:Stringy}, \cite{RuanChen:Coho},
\cite{AdemRuan:TwistedK}). 
In particular, it has become clear
that the traditional definition (originally due to Satake \cite{IS:1956})
is cumbersome and sometimes misleading. Orbifolds should not
be seen simply as manifolds with mild singularities (i.e.
a very nice kind of stratified space). Instead one wants to take
into account the extra structure provided by the stabilizer groups
attached to each point. 

From the point of view of algebraic geometry this is already known.
Orbifolds are the smooth version of Deligne-Mumford stacks.
However stacks are not well known outside of algebraic geometry,
so this point of view has not been popular among differential
geometers or topologists. When one wants to take advantage of the
extra structure of an orbifold---for example, to do ``stringy
geometry''---one must, one way or another, look at the
stack structure on the orbifold. Lupercio and Uribe have pointed
this out in \cite{LupercioUribe:TwistedK}, although they do not develop it.

More concretely, stacks are represented by groupoids, and it is
in this language that many constructions are easier than in the
traditional orbifold language. In particular the correct
notion of a map of orbifolds corresponds to a Hilsum-Skandalis
morphism of groupoids. 
Such a map induces a continuous map on the underlying topological
space, but considering only this continuous map 
loses information, that is, the functor 
$\mathbf{Orb} \ra \mathbf{Top}$ is not faithful. 
This is the precise sense in which an orbifold has ``extra structure.''
This was only observed recently in the orbifold literature, for example by Ruan
\cite{Ruan:Stringy}, who noted that to pull back an ``orbi-vector bundle,''
one needs more than a map on the underlying topological space,
and by Moerdijk and Pronk \cite{MoerdijkPronk:OrbifoldsSheavesGroupoids},
who observed the same thing for sheaves.
Also see \cite{Henriques:Borel} for concrete examples showing the need
for more data.

This paper uses the groupoid language to (partially) answer 
a question about the presentation of orbifolds. Orbifolds often arise as
a quotient of a smooth manifold $P$ by the action of 
a compact Lie group $K$. Conversely, suppose we require
that the stabilizer groups of an orbifold $\MM$ to act 
effectively on the corresponding charts. Call such
an orbifold \textit{effective}. (Some sources refer to these
as \textit{reduced} orbifolds, but this conflicts with the
standard use of the word ``reduced'' in algebraic geometry.)
Then any effective $\MM$ can be \textit{presented}
\footnote{This is very different from asking whether an orbifold
can be obtained as a quotient by the action of a finite group,
or more generally a discrete group. 
So this question has nothing to do with the distinction between
``good'' and ``bad'' orbifolds.}
as $\MM \iso P/K$ for some $P,K$. 

If one thinks of an orbifold as simply a mildly singular stratified
space, with the main function of the stabilizer groups being
to provide interesting singularities, then the effectivity assumption
is reasonable. However when one wants to utilize the extra structure
associated to the stabilizer groups, this restriction becomes
artificial. In addition, certain useful constructions on
orbifolds (such as Kawasaki's $\tilde{M}$, \cite{TK:1978}, \cite{TK:1979}) 
take one out of the effective situation.

Even the extreme case where the actions of the stabilizer groups
are all trivial turns out to be interesting; the presentation problem here 
is the topological analogue of the Brauer Conjecture 
(see \cite{EdidinHassettetal:Brauer}) in algebraic
geometry. We will see that the appropriate language for this case
uses gerbes.

Hence a natural question arises, can a noneffective orbifold
always be presented? The answer is still unknown, but we prove
results (Thm.~\ref{thm:presentbyPI}, Cor.~\ref{cor:trivialcenter} 
and Thm.~\ref{thm:presentgerbe} below) 
showing that in many cases the answer is still yes.
The extension of our method to answer the general question will require 
understanding equivariant gerbes. 

\subsection{Usual Definition of an Orbifold}\label{subsec:usual}

For completeness, we recall the usual definition of an 
orbifold. We adapt the following from Ruan and Chen \cite{RuanChen:Coho}.
(See also Moerdijk and Pronk \cite{MoerdijkPronk:OrbifoldsSheavesGroupoids}.)
We note that in Section~\ref{subsec:orbgroupoid} we will 
see an alternate characterization
of orbifolds, one that we think of as a superior definition.
Let $M$ be a paracompact Hausdorff space.

\begin{Def}\label{def:orbifold}
  An \textbf{orbifold chart} on $M$ is given by a connected open subset
  $\tilde{U} \subset \R^{n}$, an action $\rho: H \ra Aut(\tilde{U})$
  of a finite group $H$ by smooth automorphisms of $\tilde{U}$,
  and a map $\phi: \tilde{U} \ra M$, such that $\phi$ is $H$-invariant
  and induces a homeomorphism from $\tilde{U}/H$ onto an open subset
  $U = \phi (\tilde{U}) \subset M$.
  
  An \textbf{embedding} 
  $\lambda:(\tilde{U},H,\rho,\phi) \ra (\tilde{U}',H',\rho',\phi')$
  between two such charts is a smooth embedding 
  $\lambda: \tilde{U} \ra \tilde{U}'$ and an injective map, also denoted
  by $\lambda$, $\lambda: H \ra H'$, with 
  $\lambda(h \cdot \tilde{x}) = \lambda (h) \cdot \lambda (\tilde{x})$,
  and $\phi (\lambda(\tilde{x})) = \phi'(\tilde{x})$.
  Furthermore, $\lambda: H \ra H'$ is required to be an isomorphism
  from the kernel of the action map $H \ra \Aut(\tilde{U})$ to the
  kernel of $H' \ra \Aut(\tilde{U}')$. 
 
  An \textbf{orbifold atlas} on $M$ is a family 
  $\mathcal{U} = \{(\tilde{U}_{\alpha},H_{\alpha},\rho_{\alpha},\phi_{\alpha})\}$
  of such charts, which cover $M$ and are compatible in the following sense:
  given any two charts
  $(\tilde{U}_{\alpha},H_{\alpha},\rho_{\alpha},\phi_{\alpha})$ and
  $(\tilde{U}_{\beta},H_{\beta},\rho_{\beta},\phi_{\beta})$, with
  $\phi_{\alpha}(\tilde{U}_{\alpha}) = U_{\alpha}$, 
  $\phi_{\beta}(\tilde{U}_{\beta})=U_{\beta }$,
  and given any $x \in U_{\alpha} \intersect U_{\beta}$,
  there exists an open neighborhood 
  $U_{\gamma} \subset U_{\alpha} \intersect U_{\beta}$
  and a chart $(\tilde{U}_{\gamma},H_{\gamma},\rho_{\gamma},\phi_{\gamma})$
  for $U_{\gamma}$, such that there are embeddings
  $(\tilde{U}_{\gamma},H_{\gamma},\rho_{\gamma},\phi_{\gamma}) \ra 
  (\tilde{U}_{\alpha},H_{\alpha},\rho_{\alpha},\phi_{\alpha})$ and
  $(\tilde{U}_{\gamma},H_{\gamma},\rho_{\gamma},\phi_{\gamma}) \ra 
  (\tilde{U}_{\beta },H_{\beta },\rho_{\beta},\phi_{\beta})$.

  Two such atlases are said to be equivalent if they have a common refinement.
  An \textbf{orbifold} (of dimension $n$) is such a space $M$ with
  an equivalence class of orbifold atlases. 
\end{Def}

Many authors consider
only \textit{effective} (or ``reduced'') orbifolds, defined as follows.
\begin{Def}\label{def:effectiveorbifold}
  An orbifold is \textit{effective} if it has an atlas 
  in which all of the group actions $\rho_{\alpha}$ are effective.
\end{Def}
That is, for effective orbifolds, we can think of each stabilizer
group as a subgroup of $\Aut(\tilde{U})$.

For a noneffective orbifold, the role of the noneffective parts
of the stabilizer groups (the kernels of the maps $\rho_{\alpha}$)
is rather subtle. In particular the above definition is not entirely
satisfactory. In section Sec. \ref{subsec:orbgroupoid} we 
will use the language of groupoids to give
an alternate characterization of an orbifold, one which agrees
with the above definition in the effective case but which we feel
is superior in the noneffective case.

One obvious way to get an orbifold is to start with a smooth
manifold and mod out by the action of a finite group. It is
simple to find examples of orbifolds which do not arise in this
way. A more general situation in which orbifolds arise is the following
easy consequence of the slice theorem for group actions:

\begin{proposition}\label{prop:LFactiongivesorbifold}
  Let $P$ be a smooth manifold with a locally free action of
  a compact Lie group $K$. (I.e. the stabilizers are all finite.)
  Then the quotient $P/K$ has a canonical orbifold structure.
\end{proposition}

When an orbifold $M$ can be realized as such a quotient
$M \iso P/K$, we say that this is a \textit{presentation} of $M$.
In the paper \cite{IS:1956} which introduced orbifolds (under the
name ``V-manifolds'') Satake proved the following result about
presentations.

\begin{proposition}\label{prop:effpres}
  Let $M$ be an effective orbifold. Then $M$ can be presented
  as $M \iso P/K$ where $P = \mathrm{Fr}_{O(n)}(M)$ is the
  orthonormal frame bundle for some Riemannian structure on $M$
  and $K = O(n)$.
\end{proposition}

The purpose of this paper is to generalize this result to a larger
class of orbifolds than the effective ones, and give partial results
toward the presentation question for all orbifolds. 

\section{Groupoids}
\label{sec:groupoids}

As mentioned in the introduction, one of the ways to
think of an orbifold is as a particularly simple kind of
smooth groupoid. We will briefly review the notion of groupoid,
and introduce the structures on groupoids which we will need in the sequel.
A good source for some of the material in this
section is the recent preprint of Moerdijk \cite{Moerdijk:OrbGrpIntro}.
In Section~\ref{subsec:orbgroupoid} we will address how
to think of orbifolds as groupoids.

A groupoid is a category where all morphisms are invertible.
We need the notion of a smooth groupoid, where the sets of
objects and of arrows are smooth manifolds, and 
the structure maps of the groupoid are smooth maps.\footnote{One can also
work in the topological category. Here we restrict ourselves 
to the smooth setting.} 
More precisely,
\begin{Def}\label{def:groupoid}
  A \textbf{smooth groupoid} $\GG$ is a pair of
  smooth manifolds $\GG_{1}, \GG_{0}$ and
  smooth maps $s,t: \GG_{1} \ra \GG_{0}$, $u: \GG_{0} \ra \GG_{1}$, 
  $i: \GG_{1} \ra \GG_{1}$, 
  $m: \GG_{1} \times_{s,\GG_{0},t} \GG_{1} \ra \GG_{1}$,
  satisfying the relations given below.
  We call $\GG_{0}$ the space of \textbf{objects} and $\GG_{1}$
  the space of \textbf{arrows}.
  We use the notation $g: x \ra y$ to denote an arrow $g \in \GG_{1}$ with 
  $s(g) = x, t(g) = y$. Denoting $m(g,h)$ by $gh$, $i(g)$ by $g^{-1}$, 
  and $u(x)$ by $1_{x}$,
  the axioms are the familiar ones for a groupoid: for appropriate 
  $x, y \in \GG_{0}$, $g,h,k \in \GG_{1}$, 
\begin{gather}\label{groupoidstructure}
  s(1_{x}) = t(1_{x}) = x, \quad s(g^{-1}) = t(g), \quad t(g^{-1}) = s(g); \\ 
  s(gh) = s(h), \quad t(gh) = t(g), \quad g(hk) = (gh)k; \\
  g 1_{x} = 1_{y} g = g, \quad  
  g^{-1} g = 1_{x}, \quad g g^{-1} = 1_{y} \quad \text{(for $g: x \ra y$).}
\end{gather}
We further require that $s$ (and hence $t$) is a submersion (in particular, 
this guarantees that $G_{1} \times_{t,G_{0},s} G_{1}$
is actually a smooth manifold). Note: we \textbf{do} require all manifolds to
be Hausdorff (and second countable). 

If $s$ (hence $t$) is a local diffeomorphism 
we call $\GG$ a \textit{smooth \'etale groupoid}. 
If the map $(s,t): \GG_{1} \ra \GG_{0} \times \GG_{0}$ is
proper, we call the groupoid \textit{proper}.
\end{Def}

Note: It is easy to show that $i$ is an involution of $\GG_{1}$; hence 
any property holds for $s$ if and only if it holds for $t$.

When we want to refer to a groupoid without
any topological or smooth structure, we will use the term 
\textit{abstract groupoid}.

\textit{Examples.} 
\begin{enumerate}
\item Given a manifold $M$, we can define the trivial groupoid, 
      also denoted $M$,
      with $M_{0} = M$ and with only identity arrows in $M_{1}$.
\item Given a Lie group $K$ we can define a corresponding groupoid
      $\BB K$ with only one object, and with $K$ as the space of arrows.
\item Given a manifold $M$ and a right action of the Lie group $K$ on $M$, 
      we can form the translation groupoid $M \rtimes K$ by
      $(M \rtimes K)_{0} = X$, $(M \rtimes K)_{1} = M \times K$, 
      $s(x,k) = xk$, $t(x,k) = x$,
      $u(x) = (x,1)$, $i(x,k) = (xk,k^{-1})$, 
      and $(x,k) \circ (y,h) = (x,kh)$.
\end{enumerate}

Note that the first example is $M \rtimes (e)$, 
and the second example is $\ast \rtimes K$, where $\ast$, here and 
in the sequel, denotes the one-point space. The notation $\BB K$
is meant to be suggestive of classifying spaces; 
see section~\ref{subsec:HSmaps} 
for a justification of this notation.

A \textit{strict homomorphism} $\phi: \GG \ra \HHH$ 
of smooth groupoids is a smooth functor, i.e. it is 
given by a pair of smooth maps $\phi_{0}: \GG_{0} \ra \HHH_{0}$
and $\phi_{1}: \GG_{1} \ra \HHH_{1}$ commuting with the structure maps of
$\GG$ and $\HHH$. We get a corresponding notion of strict
isomorphism. However we will use a different notion of morphism,
explained in Sec.~\ref{subsec:HSmaps} below.

Given any smooth groupoid $\GG$, one can
form its topological, or coarse, quotient 
$\GG_{\mathrm{top}} = \GG_{0}/\GG$ 
by identifying any two objects which have an arrow between them.
However this is not a very fine invariant of $\GG$ and is often not
a nice topological space.
For example, any groupoid that is categorically connected (i.e.,
any two objects can be joined by an arrow) has
a topological quotient that is just a single point.
Note that if we take the topological 
quotient of a translation groupoid $X \rtimes K$, 
we get the topological quotient:
$(X \rtimes K)_{\mathrm{top}} \iso X/K$. 
We will think of groupoids as models for their quotient spaces
which are finer than the na\"{\i}ve topological quotient. 
The precise way in which we use a groupoid as a model for
its quotient comes from the notion of map which we will define
in Sec.~\ref{subsec:HSmaps}.

\subsection{Structures on Groupoids}\label{subsec:structures}

Most notions that can be defined for a group or for a manifold
can be defined for groupoids. We will need the following notions.

Let $\GG$ be a groupoid. A \textit{right action} of $\GG$ on a space
$X$ consists of two smooth maps $\pi: X \ra \GG_{0}$ (often called 
the ``moment map''; we will call it the ``base map'' to avoid
confusion with symplectic geometry usage)
and 
\begin{align*}
  m: X \times_{\pi,\GG_{0},t} \GG_{1} = \{ (x,g) \: | \: \pi(x) = t(g) \} 
           &\ra X \\
     (x,g) &\mapsto xg    
\end{align*}
(the action map) such that 
\[
  (xg)h = x(gh), \quad x1 = x, \quad \pi(xg) = s(g).
\]

Given a right action of $\GG$ on $X$ we can form the semidirect
product, or translation, groupoid, generalizing the third example
from the previous section: we define
\[
  (X \rtimes \GG)_{0} = X, \qquad 
    (X \rtimes \GG)_{1} = X \times_{\pi,\GG_{0},t} \GG_{1} 
\]
and $s(x,g) = xg$, $t(x,g) = x$, $(x,g)(xg,h) = (x,gh)$,
$i(x,g) = (xg,g^{-1})$.
Note that $\pi$ extends as a covariant functor $\pi: X \rtimes \GG \ra \GG$
by defining $\pi(x,g) = g$.
Unless otherwise specified all actions will be from the right,
so by a \textit{$\GG$-space} we will mean a space with a 
given right $\GG$-action.

A left action of $\GG$ is a right action of the opposite groupoid
$\GG^{\mathrm{op}}$. We will use this a good deal so we will be explicit.
We have a base map $\pi:X \ra \GG_{0}$ and an action map
\begin{align*}
  m: \GGl \times_{s,\GGo,\pi} X = \{ (g,x) \: | \: s(g) = \pi(x) \} 
           &\ra X \\
     (g,x) &\mapsto gx    
\end{align*}
(the action map) such that 
\[
  g(hx) = (gh)x, \quad 1x = x, \quad \pi(gx) = t(g).
\]

Again we can form a semidirect product:
\[
  (\GG \ltimes X)_{0} = X, \qquad 
    (\GG \ltimes X)_{1} = \GGl \times_{s,\GG_{0},\pi} X 
\]
and $s(g,x) = x$, $t(g,x) = gx$, $(g,hx)(h,x) = (gh,x)$,
$i(g,x) = (g^{-1},gx)$.
With these definitions $\pi$ again extends as a covariant 
functor $\pi: \GG \ltimes X \ra \GG$ by defining $\pi(g,x) = g$.

A $\GG$-sheaf of sets (or $\GG$-equivariant sheaf, or just sheaf) over
$\GG$ is a  $\GG$-space $E$ such that the base map 
$\pi: E \ra \GG_{0}$ is \'etale (i.e. a local diffeomorphism, not necessarily 
surjective). Sheaves of groups, covering spaces, and local systems
(i.e. locally constant sheaves) are defined similarly. 

We will need the following definition to describe maps of
groupoids in the next section.
\begin{Def}
  Given a manifold $B$, a \textbf{(right) $\GG$-bundle} over $B$ 
  is a (right) $\GG$-space
  $E$ and a smooth $\GG$-invariant map $p: E \ra B$ (so 
  $p(xg) = p(x)$). It is \textbf{principal} if $p$ is a 
  surjective submersion and 
  \[
    E \times_{\GG_{0}} \GG_{1} \ra E \times_{B} E, \quad (e,g) \mapsto (e,eg)
  \]
  is a diffeomorphism. 
\end{Def}
This reduces to the usual notion of a principal bundle in the
case of a Lie group, i.e. when $\GG_{0}$ is a point.

In section \ref{sec:pres}, to make certain reductions 
in the presentation problem,
we will need the notion of a group action on a groupoid
and the resulting semidirect product. (One can in fact
define the notion of a groupoid acting on a groupoid, 
but we will not need this notion.)

Let $\GG$ be a smooth groupoid and let $K$ be a Lie group. 
A \textit{right action} of $K$ on $\GG$ is a right action
(in the usual sense)
of $K$ on $\GGo$ and a right action of $K$ on $\GGl$ which
respect the structure maps of $\GG$, namely,
\[
  u(x) \cdot k = u(x \cdot k), \quad s(g \cdot k) = s(g) \cdot k,
   \quad t(g \cdot k) = t(g) \cdot k, \quad
   (gh) \cdot k = (g \cdot k) (h \cdot k).
\]
Note that every element of $K$ acts by an automorphism of 
the groupoid $\GG$.

(Note: this is a very strict notion of a group action on a
groupoid; one could define looser ``2-categorical'' 
notions but we will not need them.)

Given $K$ acting on $\GG$ on the right, we define the semidirect
product groupoid by
\[
  (\GG \rtimes K)_{0} = \GGo, \quad
    (\GG \rtimes K)_{1} = \GGl \times K
\]
with $s(g,k) = s(g) \cdot k$, $t(g,k) = t(g)$, and
$(g,k)(g',k') = (g (g' \cdot k^{-1}),kk')$.
This reduces to the usual semidirect product of groups
when $\GG$ is a group, i.e. when it has just one object;
and it reduces to the previous definition of semidirect
product when $\GG$ is a space, i.e. when it has only
identity arrows.

We have the following useful lemma about
semidirect products:
\begin{lemma}\label{lem:semidirectbytwogroups}
  Let $\GG$ be a smooth groupoid and let $K,L$ be
  two Lie groups acting on $\GG$ on the right. Assume that
  the actions commute. Then there is a natural strict isomorphism
  \[
    (\GG \rtimes K) \rtimes L \iso \GG \rtimes (K \times L). 
  \]
\end{lemma}

\subsection{Stabilizers}\label{subsec:stabilizers}

We will be particularly interested in the following $\GG$-spaces.
The \textit{stabilizer space} $S_{\GG}$ of a groupoid $\GG$ 
(denoted simply by $S$ when it does not cause confusion) is
defined to be the set of arrows with the same source and target:
\[
  S_{\GG} = \{ g \in \GG_{1} \: | \: s(g) = t(g) \}.
\]
(In categorical language, these are the automorphisms.)
There is a right action of $\GG$ by conjugation: given $g \in \GG_{1}$,
$h \in S_{\GG}$, then we define $h \cdot g = g^{-1}hg$.
Hence $S_{\GG}$ is a $\GG$-space.
Given $x \in \GGo$, let $S_{x} = \{g:x \ra x \}$ be the space of
stabilizer arrows at $x$; if $\GG$ is \'etale then $S_{x}$ is discrete.

Let $\GG$ be a smooth \'etale groupoid. An arrow $g: x \ra y$ induces
a germ of a diffeomorphism $\phi_{g}$ from $x$ to $y$. 
For, choose an open neighborhood $U$ of $g$ such that
$s|_{U}: U \ra s(U) \subset \GG_{0}$ is a diffeomorphism.
Then we have a smooth map $t \circ (s|_{U})^{-1}: s(U) \ra t(U)$,
such that $x \mapsto y$.
Different choices of $U$ give the same germ. This defines a 
strict homomorphism $\phi:\GG \ra \mathrm{Diff}(\GGo)$,
the latter being the groupoid of germs of diffeomorphisms of $\GGo$.

In particular, given a stabilizer arrow $(g:x \ra x) \in S$, consider the
germ $\phi_{g}$, which is the germ of a diffeomorphism fixing $x$.
The space of \textit{ineffective stabilizers} $S^{0} = S^{0}_{\GG}$ 
is defined as the subspace of $S$ consisting of the arrows 
which induce the trivial germ (the germ of the identity map)
on $\GG_{0}$. 
Let $S^{0}_{x} = S^{0}_{\GG} \intersect S_{x}$ be the set of 
ineffective stabilizers at $x$.

\begin{proposition}\label{prop:trivialstabilizers}
  Let $\GG$ be an \'etale groupoid. Then $S^{0}_{\GG}$ is 
  an open and closed $\GG$-equivariant subspace of $\GGl$.
  Hence $S^{0}_{\GG}$ is a $\GG$-sheaf. 
  If $\GG$ is proper, then $S^{0}_{\GG}$ is a local system.
\end{proposition}
\begin{proof}
  First note that an arrow $g \in S_{x}$ is in $S^{0}_{\GG}$ if
  and only if there is a neighborhood $U \subset \GGl$ of $g$ 
  such that $s|_{U} = t|_{U}$. Given such a $U$,
  every $h \in U$ is also in $S^{0}_{\GG}$, so $S^{0}_{\GG}$ is
  open. Also, $S^{0}_{\GG}$ is the inverse image under $\phi$ of the
  identity arrows in $\mathrm{Diff}(\GGo)$, so it is closed as well.
  
  Now we show that $S^{0}_{\GG}$ is $\GG$-equivariant.
  Given an arrow $g \in S^{0}_{x}$, let $U$ be such that $s|_{U} = t|_{U}$ 
  and there is a section $a: V = s(U) \ra U$ of $s=t$. 
  Let $h:y \ra x$ and, making $V$ smaller if necessary, let $b: V \ra \GGl$
  be a section of $t$ such that $b(x) = h$.
  Let $c:V \ra \GGl$ be defined by $c(z) = b(z)^{-1} a(z) b(z)$.
  It is easy to see that $s(c(z)) = t(c(z))$, that $c(x) = h^{-1}gh$,
  and that $c(V)$ is open in $\GGl$. Hence the conjugate $h^{-1}gh$
  is also in $S^{0}_{\GG}$.
      
  Now suppose that $\GG$ is proper. We claim that $s=t:S^{0}_{\GG} \ra \GGl$
  is proper. For, let $K \subset \GGo$ be compact. Then 
  \[
    s^{-1}(K) \intersect S^{0}_{\GG} 
       = (s,t)^{-1}(K \times K) \intersect S^{0}_{\GG}
  \]
  is the intersection of compact and a closed set, hence is compact.
  Now we claim that $s: S^{0}_{\GG} \ra \GGo$ is a covering map,
  which will complete the proof.
  This follows from the following general topology lemma.
\end{proof}

\begin{lemma}\label{lem:properetale}
  Let $p: E \ra X$ be a proper \'etale map. Assume that   
  $X$ is locally compact and $E$ is Hausdorff. Then 
  $p$ is a covering map.
\end{lemma}
\begin{proof}
  Let $x \in X$. Since $p$ is proper,
  $p^{-1}(x) = \{e_{1},\dots,e_{n}\}$ is finite,  
  and each $e_{k}$ has a neighborhood $U_{k}$ 
  such that $p|_{U_{k}}$ is a diffeomorphism onto its image.
  Since $E$ is Hausdorff, we can assume that the $U_{k}$ are
  disjoint. Let $U = \coprod U_{k}$ and let $V = \intersect p(U_{k})$.
  Since $X$ is locally compact, there is a neighborhood $W$ of
  $x$ with $\bar{W} \subset V$ compact. Let 
  \[
    A = p^{-1}(\bar{W}) \setminus U \subset E.
  \]
  which is compact, and let $B = p(A) \subset X$, also compact, hence
  closed. Since $A \intersect p^{-1}(x) = \emptyset$, 
  $x \not\in B$, so there is a neighborhood $W'$ of $x$
  with $W' \subset \bar{W}$ and $W' \intersect B = \emptyset$.
  Then $p^{-1}(W') \subset U$, so $p^{-1}(W') = \coprod U'_{k}$
  for some open subsets $U'_{k} \subset U_{k}$. Hence $W'$ is
  evenly covered, and $p$ is a covering map.
\end{proof}

If $S^{0}_{\GG}$ consists only of identity arrows, we say $\GG$
is \textit{effective}. In this paper we will be
most concerned with noneffective groupoids.

Given any groupoid $\GG$ we can quotient out by the
ineffective stabilizers. We define 
\[
  (\Geff)_{0} = \GG_{0}
   \quad \text{and} \quad 
   (\Geff)_{1} = \GG_{1}/S^{0}_{\GG}.
\]
Since $S^{0}_{\GG}$ is $\GG$-equivariant and acts trivially on $\GG_{0}$,
$\Geff$ is a smooth effective groupoid.
If $\GG$ is \'etale (resp. proper) then so is $\Geff$.
There is an evident strict homomorphism $p: \GG \ra \Geff$;
two elements $g,h:x \ra y \in \GGl$ have $p(g)=p(h)$ if 
they induce the same germ of a diffeomorphism from $x$ to $y$.

\subsection{Maps of Groupoids}
\label{subsec:HSmaps}

We have already defined the notion of a strict homomorphism of groupoids.
However, we will more often use a different notion of morphism of groupoids.
Since we are only interested in $\GG$ as a model for its
quotient, we need to be looser about what a morphism can be.
The appropriate notion of map is that of Hilsum and Skandalis.
We follow Crainic \cite{Crainic:CyclicEtale} here. 
The key is that we are (roughly) allowed to ``replace'' the source groupoid
by an ``equivalent'' object.
Here the equivalent object is a principal bundle 
over the space of objects of the source.

\begin{Def}\label{def:HSmor}
  Let $\GG$ and $\HHH$ be smooth groupoids. A (Hilsum-Skandalis) 
  morphism from $\GG$ to $\HHH$ consists of a manifold $P$, 
  maps $s_{P}: P \ra \GG_{0}$, $t_{P}: P \ra \HHH_{0}$, 
  a left action of $\GG$ on $P$
  with the base map $s_{P}$, a right action of $\HHH$ on
  $P$ with base map $t_{P}$, such that:
\begin{enumerate}
\item $s_{P}$ is $\HHH$-invariant, $t_{P}$ is $\GG$-invariant;
\item the actions of $\GG$ and $\HHH$ on $P$ are compatible: given $p \in P$,
      $(gp)h = g(ph)$;
\item $s_{P}: P \ra \GG_{0}$, as an $\HHH$-bundle with base map $t_{P}$,
      is principal. 
\end{enumerate}
We will often abuse notation and denote such a morphism simply by $P$.
\end{Def}

When we use the word ``morphism'' or ``map'' without 
qualification, we will always mean ``Hilsum-Skandalis morphism.''

The composition of two morphisms $P: \GG \ra \HHH$ and
$Q: \HHH \ra \KK$ is defined by dividing out $P \times_{\HHH_{0}} Q$
by the action of $\HHH$: $(p,q)h = (ph,h^{-1}q)$, and taking
the obvious actions of $\GG$ and $\KK$. The identity
map of $\GG$ is represented by the diagram
$\GGo \la \GGl \ra \GGo$ where the maps are the target and source
maps and the actions are the left and right multiplication.
A map $\GGo \la P \ra \HHHo$ is invertible (modulo the remark below) 
if both $s_{P}$ and $t_{P}$ are principal bundles.

\textit{Remark}. The reader may note that composition is
not associative on the nose but rather up to a natural notion
of equivalence of maps (and similarly for identity maps and inverses). 
In section \ref{sec:stacks}
we will note some of the subtleties this causes, especially  
relating to gluing. For the most part these issues will
not concern us. However as one reflection of this subtlety
we prefer to call a map with $s_{P},t_{P}$ both principal
an \textit{equivalence} rather than an isomorphism.

Given a strict homomorphism $\phi: \GG \ra \HHH$, there is
an associated HS morphism given by
\[
  \xymatrix{ \GGo & \GGo \times_{\phi, \HHHo, t} \HHHl 
                      \ar[l]_-{\mathrm{pr}_{1}} 
                      \ar[r]^-{s \circ \mathrm{pr}_{2}} 
                  & \HHHo
  }
\]
where the actions of $\GG$ and $\HHH$ on $\GGo \times_{\HHHo} \HHHl$
are given by $g \cdot (x,h) \cdot h' = (t(g),\phi(g) h h')$.
It is easy to check that this defines an HS morphism.
The resulting HS morphism will be an equivalence exactly when
the functor giving rise to it is an \textit{essential equivalence},
defined as follows.
\begin{DefProp}\label{def:essequiv}
  A strict homomorphism $\phi: \GG \ra \HHH$ is an
  \textbf{essential equivalence} if one, hence both, of the following
  equivalent conditions are satisfied:
  \begin{enumerate}
  \item the resulting HS map is an HS equivalence;
  \item \label{esseq}
        \begin{enumerate}
        \item \label{esssurj} $\phi$ is \textbf{essentially surjective:}
              the map 
              \[
                 s \circ \mathrm{pr}_{2}:\GGo \times_{\phi,\HHHo,t} 
                                                       \HHHl \ra \HHHo
              \]
              is a surjective submersion; 
        \item \label{fullyfaithful} $\phi$ is \textbf{fully faithful:} 
              the square
              \[
                \xymatrix{\GGl \ar[r]^-{\phi} \ar[d]_-{(s,t)}
                          & \HHHl \ar[d]^-{(s,t)} \\ 
                          \GGo \times \GGo \ar[r]^-{\phi \times \phi}    
                          & \HHHo \times \HHHo
                }
              \]
              is a fiber product.
        \end{enumerate}
\end{enumerate}
\end{DefProp}
It is easy to check the equivalence of these conditions.

Note that the second condition implies that the underlying map of
categories (without any topological or smooth structure) is
an equivalence.

\textit{Remark.} One can arrive at the notion of morphism of groupoids
by inverting the essential equivalences (in the sense of localization
of categories). However we prefer the HS definition, partly
because it leads more easily to the ``stacky'' nature
of groupoids, see Section~\ref{sec:stacks}. 
(See also \cite{Pronk:Bicategories}.)

\textit{Example.} 
We can let the source and the target of an HS morphism be  
translation groupoids. 
Let $X$ be a $K$-space, $Y$ be an $L$-space, for Lie groups $K,L$.

\begin{lemma}\label{HSintranslationcase}
  A morphism of translation groupoids $X \rtimes K \ra Y \rtimes L$
  is equivalent to the data of a $K$-equivariant 
  principal $L$-bundle $s_{P}: P \ra X$ over the $K$-space $X$
  and a $K$-invariant and $L$-equivariant map $t_{P}: P \ra Y$.
\end{lemma}
(Note: since $K$ acts on the right of $X$ and on the left of  
$P$, $s_{P}$ being $K$-equivariant means 
$s_{P}(k \cdot p) = s_{P}(p) \cdot k^{-1}$.

The simple proof is left to the reader. 
Note that if we understand $X$ and $Y$ to
be given the trivial actions of $L$ and $K$ respectively, 
we can describe both $s_{P}$ and $t_{P}$ as 
$K^{\mathrm{op}} \times L$-equivariant maps. 

We can further specialize by setting $Y$ to be a point, $Y=\ast$, 
and $K$ to be the trivial group, $K=(e)$. 
Then a map from $X = X \rtimes (e)$ to $\BB L = \ast \rtimes L$ is simply a
principal $L$-bundle over $X$. in other words, $\BB L$ classifies
$L$-bundles, as the notation suggests. 
(Note that this fits with the topological
definition of the classifying space, which is $BG = EG/G$.
$EG$ is contractible, hence ``morally'' a point; topologists use it
instead of a point to get a ``good'' quotient. Here we actually use
a point for the space of objects, 
but use an alternate category of ``spaces'' to take the
quotient.) We will discuss this example further in Section~\ref{sec:stacks}.

Many (but not all) properties and invariants of a groupoid
are preserved by HS equivalence. For example, we have the following,
saying that topological quotients and stabilizer groups are
preserved.
\begin{proposition}\label{prop:HSeqfeatures}
  Let $(P,s_{P},t_{P})$ be an HS map from $\GG$ to $\HHH$.
  Then $P$ induces a continuous map on the underlying topological
  quotients, $f_{P}:\Gtop \ra \Htop$. Given $x \in \GGo$ and
  $y \in \HHHo$ such that $f_{P}([x]) = [y]$, $P$ induces a group
  homomorphism $\psi_{P,x,y}: S_{x} \ra S_{y}$, well defined up to
  conjugacy in the target.

  If $P$ is an equivalence, then $f_{P}$ is a homeomorphism
  and each $\psi_{P,x,y}$ is an isomorphism.
\end{proposition}
\begin{proof}
  The map $f_{P}:\Gtop \ra \Htop$ is defined by $f_{P}([s_{P}(p)]) = [t_{P}(p)]$.
  Since $s_{P}$ is surjective, this defines $f_{P}$ for all $[x] \in \Gtop$.
  To check that it is well-defined, let $[s_{P}(p)]=[s_{P}(q)]$;
  then there is a $g:s_{P}(p) \ra s_{P}(q)$ in $\GG$, so $s_{P}(g \cdot p) = s_{P}(q)$.
  Since the $\HHH$-action on $P$ is principal, there is $h \in \HHHl$ with
  $g \cdot p \cdot h = q$ and $h: t_{P}(q) \ra t_{P}(g \cdot p)$.
  Since $t_{P}$ is $\GG$-invariant, we have $h:t_{P}(q) \ra t_{P}(p)$,
  so $[t_{P}(p)]=[t_{P}(q)]$. It is also straightforward to check that
  $f_{P}$ is continuous. If $P$ is an equivalence, then the construction
  works in both directions, and it simple to check that the resulting
  maps are inverses.

  Given $x \in \GGo$, $y \in \HHHo$ with $f_{P}([x]) = [y]$, we 
  define $\psi_{P,x,y}: S_{x} \ra S_{y}$ as follows. Let 
  $p \in P$ with $s_{P}(p) = x$. We know there is 
  $h_{0}:t_{P}(p) \ra y$; fix such an $h_{0}$.
  Given $g \in S_{x}$, there is a unique $h \in HHHl$ with $ph = gp$; 
  here $h: t_{P}(p) = t_{P}(gp) = t_{P}(ph) \ra t_{P}(p)$. 
  Hence we define $\psi_{P,x,y}(g) = h_{0} h h_{0}^{-1} \in S_{y}$.
  It is simple to check that changing $p$ or $h_{0}$
  changes the map $\psi_{P,x,y}$ by conjugating in the target.
  If $P$ is an equivalence, as above, the construction works
  in both directions to give inverse isomorphisms.
\end{proof}

\subsection{Orbifold Groupoids}\label{subsec:orbgroupoid}

Moerdijk and Pronk (\cite{MoerdijkPronk:OrbifoldsSheavesGroupoids}, 
\cite{Moerdijk:OrbGrpIntro})
show that to an effective orbifold with a given atlas, 
one can naturally associate
a smooth proper effective \'etale groupoid. Equivalent
atlases give rise to equivalent groupoids.
Every such groupoid arises from an effective orbifold.
Further, Ruan's good maps (\cite{Ruan:Stringy}) of orbifolds correspond to
HS maps of groupoids. 

Hence instead of the traditional definition of orbifold,
we will henceforth consider orbifolds as groupoids.
More precisely, following Moerdijk, 
we will refer to a smooth proper \'etale groupoid as an
\textit{orbifold groupoid}, and refer to a smooth proper effective \'etale
groupoid as an \textit{effective orbifold groupoid}. 
(Actually Moerdijk allows the groupoid to be an arbitrary proper foliation
groupoid, which is always equivalent to an \'etale groupoid.)
Similarly, a map of orbifold groupoids (our replacement for the notion of
a good map of orbifolds) will be a Hilsum-Skandalis map.

We remark that this gives an elegant answer to the question of
exactly what conditions one should put on the ineffective parts
of the stabilizers, which we alluded to in Section \ref{subsec:usual}.
One can translate the groupoid statement to one involving explicit charts,
analogous to Definition \ref{def:orbifold}, but it is not particularly
enlightening. One of our main themes in this paper is to stay
within the groupoid language and demonstrate its usefulness.

Henceforth we will freely use the word ``orbifold'' to mean
``orbifold groupoid.''

If we have presented orbifolds
$M = X/K$, $N = Y/L$, then by Lemma~\ref{HSintranslationcase},
a map between the orbifolds is given by 
a $K$-equivariant principal $L$-bundle $s_{P}: P \ra X$
and a $K$-invariant, $L$-equivariant map $t_{P}: P \ra Y$.

\section{Groupoids as Stacks}
\label{sec:stacks}

We want to point out some of the subtleties inherent in 
the definition of maps of orbifolds given in Section~\ref{subsec:HSmaps},
and briefly relate orbifolds to stacks. This section
can be skipped if the reader is only interested in
the presentation problem. However it does serve as motivation
for the next section on gerbes.

First, given two HS maps from $\GG$ to $\HHH$, one can have 
an ``isomorphism of HS maps'' between them. It turns out that
we do not want to simply identify two such isomorphic maps, as
we would then lose the ability to glue maps. This leads to the
idea of thinking of a groupoid as representing a \textit{stack}.
This point of view is very natural to algebraic geometers,
but less so to differential geometers and topologists. However
it should underlie all of the new ``stringy'' constructions on
orbifolds. (See \cite{LupercioUribe:TwistedK} for more comments on this.)

\begin{Def}\label{def:HS2iso}
  Given two groupoids $\GG$ and $\HHH$ and two HS maps 
  $(P,s_{P},t_{P}): \GG \ra \HHH$ and $(Q,s_{Q},t_{Q}): \GG \ra \HHH$, 
  an \textbf{isomorphism} (or \textbf{2-isomorphism} for definiteness)
  from $(P,s_{P},t_{P})$ to $(Q,s_{Q},t_{Q})$ is a 
  $\GG \times \HHH$-equivariant diffeomorphism $\alpha: P \ra Q$
  such that $s_{Q} \alpha = s_{P}$ and $t_{Q} \alpha = t_{P}$.
\end{Def}
Hence for every $\GG, \HHH$ there is a category 
(in fact, an abstract groupoid) 
$\HS(\GG,\HHH)$ of HS maps from $\GG$ to $\HHH$. 

One would like to say this makes (Groupoids, HS maps) into
a 2-category, but that is not quite true. Since the
definition of composition of HS maps involves taking a pullback,
composition of morphisms is not strictly associative, but only
associative up to a coherent isomorphism. 
Hence (Groupoids, HS maps) form what is called a \textit{bicategory}. 
However the lack of strict associativity of composition will not bother
us too much, so it is not too wrong to think of this as a 2-category.

It is tempting to identify any two isomorphic HS maps to get
a set of equivalence classes of maps. This does form a 
category, and often this is good enough. However there is one
fundamental operation that fails when one does this, namely 
gluing maps together. A simple example shows the basic problem.

As in the previous section, let 
$M$ be a manifold and let $K$ be a group and consider the category
of HS maps from the trivial groupoid $M$ to the classifying groupoid
$\BB K$. We saw that $HS(M,\BB K)$ is the category of $K$-principal bundles
on $M$. Suppose that $M = \bigcup_{\alpha} U_{a}$ is an open cover
of $M$. Then a $K$-bundle over $M$ is determined by the
data of a bundle on each $U_{\alpha}$ and specific gluing 
isomorphisms on each $U_{\alpha} \intersect U_{\beta}$,
such that the cocycle condition is satisfied.

If we simply identified any two isomorphic HS maps,
then we would lose these specific gluing isomorphisms,
retaining only the knowledge that the bundles happen to be 
isomorphic. If the cover is fine enough, that is no information
at all since all bundles are locally trivial. Also, there
is no way to talk about compatibility of gluings.

Hence HS maps glue in a different way from ordinary maps of
manifolds. We have the following setup.

Let $M$ be a manifold (considered as a trivial groupoid)
and $\GG$ be an arbitrary groupoid. Let 
$\UU = \{U_{\alpha}\}$ be an open
cover of $M$, and denote $U_{\alpha} \intersect U_{\beta}$ 
by $U_{\alpha \beta}$. 
\begin{Def}\label{def:descentdata}
  A collection of HS maps 
  \[
    \psi_{\alpha}: U_{\alpha} \ra \GG
  \]
  and 2-isomorphisms of HS maps 
  \[
    \chi_{\alpha \beta}: \psi_{\alpha}|_{U_{\alpha \beta}} \ra 
       \psi_{\beta}|_{U_{\alpha \beta}} 
  \]
  is called \textbf{descent data} if the $\chi_{\alpha \beta}$
  satisfy the cocycle condition
  \[
    \chi_{\alpha \beta} \chi_{\beta \gamma} \chi_{\gamma \alpha} = 1.
  \]
  We will only consider decent data that are \textbf{normalized}, i.e. 
  with 
  \[
    \chi_{\beta \alpha} = (\chi_{\alpha \beta})^{-1}
     \quad \text{and} \quad \chi_{\alpha \alpha} = 1.
  \]
  Descent data with respect to a given open cover form
  a groupoid, denoted $\Desc(\UU,\GG)$. 
  An \textbf{isomorphism} of descent data from
  $(\psi,\chi)$ to $(\psi',\chi')$ is given by 
  \[
    \eta_{\alpha}: \psi_{\alpha} \ra \psi'_{\alpha}
  \]  
  such that 
  \[
    \chi'_{\alpha \beta} \circ \eta_{\alpha}|_{U_{\alpha \beta}} 
      = \eta_{\alpha}|_{U_{\alpha \beta}} \circ \chi_{\alpha \beta}
  \]
\end{Def}

The proof of the following proposition exactly parallels the
proof of the corresponding fact for principal bundles.
\begin{proposition}\label{prop:HSmapsarestacky}
  Let $M$ be a manifold (considered as a trivial groupoid)
  and let $\GG$ be a groupoid. Let $\UU = \{U_{\alpha}\}$ be an open
  cover of $M$. For every HS map $\psi: M \ra \GG$,
  restriction to the open cover induces an equivalence
  between the categories
  \[
    HS(M,\GG) \iso \Desc(\UU,\GG).
  \]
\end{proposition}

In the algebraic geometry language, 
this says that a groupoid $\GG$ naturally induces
a \textit{stack} on the category of manifolds.
We present some further explanation of this idea
and clarification of some confusing issues
in a related paper \cite{Metzler:Stacks}.

\section{Gerbes}
\label{sec:gerbes}

We will be particularly interested in a very simple kind of orbifold groupoid $\GG$,
namely one whose corresponding effective orbifold $\Geff$ is equivalent to 
a manifold. 

\begin{Def}\label{def:PIgroupoid}
  A smooth \'etale groupoid $\GG$ is \textbf{purely ineffective (PI)}  
  if $S^{0}_{\GG} = S_{\GG}$, i.e. if every stabilizer arrow gives the
  germ of the identity on $\GGo$.
\end{Def}

\begin{proposition}\label{prop:PIgivesetale}
  Let $\GG$ be a purely ineffective orbifold groupoid. Then the quotient map
  $p: \GGo \ra \Gtop$ is \'etale and the quotient $\Gtop$
  is a smooth manifold. Moreover the induced map from $\Geff$ to $\Gtop$
  is an equivalence.
\end{proposition}
\begin{proof}
  First note that the quotient map is open, since given an open
  set $U \subset \GGo$, $p^{-1}p(U) = s t^{-1}(U)$ is open since 
  $s$ is an open map.

  Further, note that $S = S^{0}$ is open in $\GGl$, by 
  Prop.~\ref{prop:trivialstabilizers}. 

  Let $x \in \GGo$. We will show that there is 
  a neighborhood $W$ of $x$ such that
  for any arrow $g:y \ra z$, $y,z \in W$, we must have $y=z$.
  Given such a $W$, $p|_{W}$ is open and injective, hence a homeomorphism
  onto its (open) image. This will prove $p$ \'etale. It will follow 
  that $\Gtop$ is a smooth manifold (it is Hausdorff because
  it is a quotient by a proper equivalence relation).

  Now we show the existence of $W$. 
  Let $U$ be a neighborhood of $x$ whose closure 
  is compact. Let $A = (s,t)^{-1}(\bar{U} \times \bar{U})$.
  Since $\GG$ is proper, $A$ is compact. Let $B = A \setminus S$
  which is also compact. 

  Suppose that there is a sequence of points $y_{i},z_{i} \in U$ with
  arrows $h_{i}: x_{i} \ra z_{i}$, such that $y_{i}, z_{i} \ra x$,
  and $h_{i} \not\in S$. We have $h_{i} \in B$ which is compact,
  so there is a convergent subsequence $h_{i_{j}} \ra h \in B$ of $h_{i}$.
  This $h$ must be an arrow with $s(h) = t(h) = x$, so $h \in S$,
  a contradiction. Hence no such sequence can exist, so there
  must be some neighborhood $W$ of the form desired; hence
  $p$ is \'etale. 
  
  To show that the induced strict homomorphism, 
  which we shall also call $p$, from
  $\Geff$ to $\Gtop$ is an equivalence, we need to verify 
  conditions (\ref{esssurj}), (\ref{fullyfaithful})  
  in Def-Prop.~\ref{def:essequiv}. Essential surjectivity is clear
  since the map on objects is a surjective \'etale map, hence certainly
  a surjective submersion. The map $p$ is also fully faithful
  since for both groupoids, there is a unique arrow from
  any object to any other object.
\end{proof}

As a corollary, we know when a groupoid is equivalent to a 
manifold. The answer is simple and standard.
\begin{cor}\label{cor:whenmanifold}
  Given a smooth proper \'etale groupoid $\GG$, there is a manifold
  $M$ whose associated groupoid is equivalent to $\GG$
  if and only if $S_{\GG}$ consists only of identity arrows.
  In that case, $\Gtop$ is a smooth manifold and is diffeomorphic to $M$.
\end{cor}
\begin{proof}
  Given a smooth \'etale groupoid with $S_{\GG}$ consisting only
  of identity arrows, Prop.~\ref{prop:PIgivesetale} implies that
  $\Gtop$ is a smooth manifold. It is easy to see that $p$
  satisfies condition~\ref{esseq} of the definition of an essential
  equivalence.

  Conversely, since stabilizer groups are preserved by equivalence,
  if $\GG$ is equivalent to a manifold, it must have no nonidentity
  stabilizers.
\end{proof}

Prop.~\ref{prop:PIgivesetale} says that a PI orbifold groupoid is not
a singular space at all, but instead is a smooth manifold with some
extra structure, given by the stabilizers of the groupoid. In fact
the extra structure is a well-known one, at least to algebraic
geometers. Such a groupoid corresponds to a \textit{gerbe} on the
underlying manifold $\Gtop$.

We recall the definition of a gerbe $F$ on a space $M$. 
See \cite{Giraud:Book}, \cite{Brylinski:LoopBook},
\cite{LMB} for more details.
A gerbe is roughly a certain kind of ``sheaf of categories.''
(More precisely, it is a certain kind of stack.)
\begin{Def}\label{def:gerbe}
  Let $M$ be a topological space. A (strict) \textbf{gerbe} $F$ on $M$
  is an assignment of an abstract groupoid $F(U)$ to every open set
  $U \subset M$ and a strict homomorphism of groupoids (a functor) 
  $r_{VU}:F(U) \ra F(V)$ for every inclusion $V \subset U$.
  Denote $r_{VU}(x) = x|_{V}$ for $x \in \Ob F(U)$,
  and for $x,y \in \Ob F(U)$, denote the set of arrows in $F(U)$ 
  from $x$ to $y$ by $F(x,y)$. This data is required to satisfy:
  \begin{enumerate}
\item \label{gerbepresheaf} 
      $r_{WV} \circ r_{VU} = r_{WU}$ for $W \subset V \subset U$.
\item \label{gerbeprestack}
      Given $U \subset M$ and $x,y \in \Ob F(U)$, 
      the assignment $V \mapsto F(x|_{V},y|_{V})$
      defines a sheaf (of sets).
\item \label{gerbestack}
      Suppose we are 
      given an open cover $\mathcal{U} = \{U_{\alpha}\}$ 
      of $U$ and \textbf{descent data}
      for $F$ over $\mathcal{U}$, i.e. $x_{\alpha} \in \Ob F(U_{\alpha})$
      and $\phi_{\alpha \beta}: x_{\alpha}|_{U_{\alpha} \intersect U_{\beta}} 
      \iso x_{\beta}|_{U_{\alpha} \intersect U_{\beta}}$ such that
      $\phi_{\alpha \beta} \phi_{\beta \gamma} = \phi_{\alpha \gamma}$
      (when restricted to 
      $U_{\alpha} \intersect U_{\beta} \intersect U_{\gamma}$).
      Then there is an $x \in \Ob F(U)$ 
      and $\psi_{\alpha}: x|_{U_{\alpha}} \iso x_{\alpha}$
      such that $\phi_{\alpha \beta} \psi_{\alpha} = \psi_{\beta}$ 
      (over $U_{\alpha} \intersect U_{\beta}$).
\item \label{gerbenonempty}
      There is an open cover $\mathcal{U} = \{U_{\alpha}\}$ of $M$
      such that each category $F(U_{\alpha})$ is nonempty.
\item \label{gerbelocaliso}
      Given $U \subset M$ and $x,y \in \Ob F(U)$, there is an open cover
      $\mathcal{V} = \{V_{\alpha}\}$ of $U$ such that for every $\alpha$,
      $x|_{V_{\alpha}} \iso y|_{V_{\alpha}}$.
\end{enumerate}
\end{Def}
This definition amalgamates a number of more general concepts.
Condition (\ref{gerbepresheaf}) says that $F$ is a presheaf of categories.
Adding condition (\ref{gerbeprestack}) says that $F$ is a ``prestack.''
Adding condition (\ref{gerbestack}) says that $F$ is in fact a stack.
One refers to this condition by saying that all descent data
are \textit{effective}, i.e. they define objects of $F(U)$ as stated.

Condition (\ref{gerbenonempty}) says that the categories $F(U)$ are
locally nonempty, and condition (\ref{gerbelocaliso}) says that
any two objects are locally isomorphic. These last two conditions
say that the stack $F$ is in fact a gerbe.

\textit{Remark.} Condition (\ref{gerbepresheaf}), that the restriction
functors commute on the nose, makes this 
a \textit{strict} gerbe. Usually one allows  
the restriction functors to commute up to a coherent natural isomorphism.
However in our case the strict condition will be satisfied, which simplifies
the notation.

Now suppose $\GG$ is a purely ineffective orbifold groupoid with topological
quotient $\Gtop$. 

We define a corresponding gerbe $F$ over $\Gtop$
as follows. Informally this is done as follows. 
Locally, over a small enough open set $U \subset \Gtop$,
\begin{equation*}
  \Ob F(U) = \Gamma(U;\GGo) 
           = \{ \text{sections of~} p:\GG_{0} \ra \Gtop \text{~over~} U\}.
\end{equation*} 
The arrows are defined locally, for two local sections 
$\sigma,\tau: U \ra \GG_{0}$, by 
\[
 \Ar F(U)(\sigma,\tau) 
   = \Gamma(U;\GGl)_{\sigma,\tau} = \{ \gamma: U \ra \GG_{1} \: | \: 
        s \circ \gamma = \sigma, t \circ \gamma = \tau \}.
\]
Informally, the arrows over $U$ are just the different ways of
gluing the representatives of $U$ together.

Formally, we first define a sheaf of categories $\hat{F}$ by 
\begin{align*}
  \Ob \hat{F}(U) &= \{ \text{sections of~} p:\GG_{0} \ra \Gtop\}, \\
  \Ar \hat{F}(U)(s,t) 
    &= \{ \gamma: U \ra \GG_{1} \: | \: 
        s \circ \gamma = \sigma, t \circ \gamma = \tau \}
\end{align*}
and then take $F$ to be the \textit{associated stack}.
This is defined as follows (following \cite{LMB}. 
An object of $F(U)$ consists of
a cover $\UU = \{U_{\alpha}\}$ of $U$, and descent data over $\UU$, i.e.
a family $x_{\alpha} \in \hat{F}(U_{\alpha}) = \Gamma(U_{\alpha};\GGo)$,
and morphisms 
\[
  f_{\alpha \beta}: x_{\beta}|_{U_{\alpha}\intersect U_{\beta}}
   \ra x_{\alpha}|_{U_{\alpha}\intersect U_{\beta}}
\]
satisfying $f_{\alpha \beta} f_{\beta \gamma} = f_{\alpha \gamma}$.  
By the definition of $\hat{F}$, 
$f_{\alpha \beta}$ is a section of $\GGl$ over $U_{\alpha} \intersect U_{\beta}$.

An arrow of $F(U)$ from $(\{U_{\alpha}\},\{x_{\alpha}\},\{f_{\alpha \beta}\})$
to $(\{U'_{\gamma}\},\{x'_{\gamma}\},\{f'_{\gamma \delta}\})$
is given by the following data. 
Let $V_{\alpha \gamma} = U_{\alpha} \intersect U'_{\gamma}$; these form
a cover of $U$ which refines both $\{U_{\alpha}\}$ and $\{U'_{\gamma}\}$. 
Then an arrow consists of a family of morphisms $\{g_{\gamma \alpha}\}$ of 
$\hat{F}(V_{\alpha \gamma})$ from
$x_{\alpha}|_{V_{\alpha \gamma}}$ to $x'_{\gamma}|_{V_{\alpha \gamma}}$,
compatible with the descent isomorphisms:
\[
   g_{\gamma \alpha} f_{\alpha \beta} = f'_{\gamma \delta} g_{\delta \beta}.
\]
(This is defined over 
$U_{\alpha} \intersect U_{\beta} \intersect U'_{\gamma} \intersect U'_{\delta}$).

Restriction of objects and arrows from $F(U)$ to $F(V)$, for $V \subset U$, is
defined by restriction of $x_{\alpha}$, etc., in the obvious way.

It is straightforward to verify that this satisfies 
conditions (\ref{gerbepresheaf})--(\ref{gerbestack}) 
of Def.~\ref{def:gerbe}. 

By Prop.~\ref{prop:PIgivesetale},
$\GG_{0} \ra \Geff = \Gtop$ always has local sections (being \'etale),
verifying condition (\ref{gerbenonempty}).
Further, the map $\GG_{1} \ra \GG_{0}$ has local sections,
and hence any two objects of $F(U)$ are locally isomorphic,
verifying condition (\ref{gerbelocaliso}).
Therefore $F$ is indeed a gerbe. Note that if there are effective
stabilizers, one will not get local sections of $\GG_{0} \ra \Gtop$,
so one does not obtain a gerbe over $\Gtop$.

Conversely, if we start with a smooth manifold $M$ and a 
certain type of gerbe $F$ on $M$ (one whose band is locally constant,
with finite isotropy groups)
one can construct a corresponding orbifold groupoid $\GG = (M,F)$.
In fact there is an equivalence between orbifolds with
only ineffective stabilizers and these types of gerbes over manifolds.
We will not need this result here.


\section{Presentations of Orbifolds}
\label{sec:pres}

In this section we will prove a number of propositions that
reduce the presentation problem for a general orbifold groupoid
to a particular problem about equivariant abelian gerbes.
We also show that the corresponding nonequivariant problem
is solvable.

First we desingularize the space using the standard procedure
of passing to the frame bundle.

Let $\GG$ be a smooth \'etale groupoid of dimension $n$.
Let $F = \mathrm{Fr}(T \GG_{0})$ be the orthonormal frame bundle
of $\GG_{0}$. Then there is a natural left action of $\GG$ on
$F$ with the standard projection $\pi: F \ra \GGo$ as the base map.
For, since $\GG$ is \'etale, given an arrow $g: x \ra y$
and a frame $f = (e_{1},\dots ,e_{n})$ with each $e_{k} \in T_{y}(\GGo)$,
we can use the induced germ of a diffeomorphism $\phi_{g}$
to act on the frame: $g \cdot f = \phi_{g}(f)$.
This is easily seen to be a smooth left action of $\GG$.
Further, this action commutes with the natural right action
of $O(n)$ on $F$.

We can therefore use the following, which is stated with
an extra group action present for later use.

\begin{lemma}\label{lem:principalbundle}
  Let $\GG$ be a smooth groupoid and let $K,L$ be Lie groups.
  Let $\pi: F \ra \GGo$ 
  be a left $\GG$-equivariant, right $L$-principal bundle 
  over $\GGo$, and assume that the  
  actions of $\GG$ and $L$ commute.
  Further, assume that: $K$ acts on $\GG$ and on $F$,
  $\pi$ is $K$ equivariant, the actions of $K$ and $L$ commute,
  and the actions of $\GG$ and $K$ satisfy
  \[
    (g \cdot f) \cdot k = (g \cdot k) \cdot (f \cdot k)
  \]
  for $g \in \GGl, f \in F, k \in K$.
  Then $K \times L$ acts on the right of $\GG \ltimes F$ and 
  the map $\pi$ induces a $K$-equivariant equivalence between
  $(\GG \ltimes F) \rtimes L$ and $\GG$.
\end{lemma}
\begin{proof}
  It is clear that $K \times L$ acts on the right of $\GG \ltimes F$,
  by $(g,f) \cdot (k,l) = (g \cdot k, (f \cdot k) \cdot l)$.
  We have 
  \[
    ((\GG \ltimes F) \rtimes L)_{0} = F, \quad
      ((\GG \ltimes F) \rtimes L)_{1} = \GGl \times_{s,\GGo,\pi} F \times L.
  \]
  We extend $\pi: F \ra \GGo$ to arrows by defining   
  $\pi(g,f,l) = g$. As usual this defines a smooth covariant functor,
  i.e. a strict homomorphism of groupoids. On objects,
  $\pi$ is a surjective submersion, since it is a principal bundle. 
  We also need to check that $\pi$ is fully faithful.
  Since $F$ is an $L$-principal bundle, we have 
  $F \times L \iso F \times_{\GGo} F$. Hence
  \begin{align*}
    \GGl \times_{\GGo} F \times L 
       &\iso \GGl \times_{\GGo} F \times_{\GGo} F \\
       &\iso \GGl \times_{\GGo \times \GGo} (F \times F) 
  \end{align*}
  which says that $\pi$ is fully faithful.

  Hence $\pi$ is an essential equivalence. It is also
  clearly $K$-equivariant. 
\end{proof}

Note that $\GG \ltimes F$ should be
thought of as a principal $L$-bundle over $\GG$. Also, we want to
think of $(\GG \ltimes F) \rtimes L$ as (a finer model of) the quotient
of $\GG \ltimes F$ by the action of $L$. So the lemma can be thought
of as simply saying that when we mod out the total space of an
$L$-bundle by $L$, we recover the base space.

Our first application of the lemma is the following.

\begin{proposition}\label{prop:framebundle}
  Let $\GG$ be an orbifold groupoid of dimension $n$.
  Put a $\GG$-invariant Riemannian metric on $\GGo$
  and let $\pi: F = \mathrm{Fr}(T \GGo) \ra \GGo$ be the orthonormal frame bundle
  of $\GGo$. Then $\GG \ltimes F$ is a purely ineffective
  groupoid with a right action of $O(n)$,
  and $\pi$ induces an $O(n)$-equivariant equivalence
  $(\GG \ltimes F) \rtimes O(n) \iso \GG$.
  Further, $S^{0}_{(\GG \ltimes F)} \iso \pi^{*}(S^{0}_{\GG})$.
\end{proposition}
\begin{proof}
  The action of $O(n)$ and the stated equivalence follow 
  directly from the lemma (with $K$ trivial). We need only
  to show that $\GG \ltimes F$ is purely ineffective.
  Given $x \in \GGo$, $f = (e_{1},\dots, e_{n})$, $e_{k} \in T_{x}\GGo$,
  and $g \in \GGl$, $g: x \ra y$, suppose that 
  $(g,f) \in S_{(\GG \ltimes F)}$, i.e. that 
  $g \cdot f = f$. Then $(\phi_{g})_{*} = 1$ on $T_{x}\GGo$.
  Now $g \in S_{\GG}(x)$ which is a finite group, so
  $\phi_{g}$ must be the identity. (Use an invariant inner product
  on a neighborhood of $x$ in $\GGo$ and the resulting exponential map.) 
  Hence $g \in S^{0}_{\GG}$, and $(g,f) \in S^{0}_{(\GG \ltimes F)}$.
  Hence 
  $S_{(\GG \ltimes F)} = S^{0}_{(\GG \ltimes F)} \iso \pi^{*}(S^{0}_{\GG})$.
\end{proof}

In particular, if $\GG$ is effective, then $\GG \ltimes F$
has no stabilizers, so is equivalent to a manifold. This
is the groupoid version of the 
proof that an effective orbifold is presentable
as the quotient of its frame bundle by $O(n)$.

Now we reduce the complexity of $\GG$ at the expense of 
making the compact group larger.
Given a purely ineffective orbifold groupoid $\GG$, the stabilizers will
not generally be abelian, and 
the local system $S^{0}_{\GG}$ may not be $\GG$-equivariantly trivial.
We will replace $\GG$ with a new groupoid $\GG'$ where $S^{0}_{\GG'}$ is
abelian and $\GG$-equivariantly trivial, using a standard trick from bundle theory.
We will assume $\GG$ is connected (i.e. $\Gtop$ is connected)
for simplicity.

\begin{proposition}\label{prop:trivialband}
  Let $\GG$ be a purely ineffective connected orbifold groupoid and
  let $K$ be a Lie group acting on $\GG$ on the right. 
  Let $T = S^{0}_{\GG}(x)$ for some $x \in \GGo$.
  Then there is a purely ineffective orbifold groupoid 
  $\GG'$ with an action of $K \times \Aut(T)$ such that
  $\GG' \rtimes \Aut(T)$ is $K$-equivariantly equivalent to $\GG$
  and the stabilizer local system $S^{0}_{\GG'}$ of $\GG'$ 
  is ($\GG'$,$K \times \Aut(T)$)-equivariantly 
  isomorphic to $\GGo' \times Z(T)$.
  In particular the stabilizer local system of $\GG'$ is abelian and 
  $K$-equivariantly trivial.
\end{proposition}
\begin{proof}
  Let 
  \[
    F = \{ (x,\phi) \in \GGo \times \Hom(T,\GG) 
                \: | \: \phi: T \iso S^{0}_{\GG}(x) \}.
  \]
  There is an obvious projection map $\pi: F \ra \GGo$.
  It is surjective since all of the stabilizer groups
  of $\GG$ are isomorphic to $T$ (since $\Gtop$ is connected and $S^{0}_{\GG}$ 
  is a local system). 
  In fact it makes $F$ into a 
  $\GG$-equivariant $\Aut(T)$-principal
  bundle, where the left action of $(g: x \ra y) \in \GG$ is
  \[
    g \cdot (x,\phi) = (y,g \phi g^{-1})
  \]
  and the right action of $\lambda \in \Aut(T)$ is 
  \begin{equation}\label{autTaction}
    (x,\phi) \cdot \lambda = (x, \phi \circ \lambda).
  \end{equation}
  These actions clearly commute.
  It is also easy to check that everything is $K$-equivariant.
  Hence by Lemma~\ref{lem:principalbundle}, $K \times \Aut(T)$
  acts on the right of $\GG \ltimes F$ and $\pi$ induces
  a $K$-equivariant equivalence between 
  $(\GG \ltimes F) \rtimes \Aut(T)$ and $\GG$.
  So we let $\GG' = \GG \ltimes F$.

  It is clear that $S^{0}_{\GG'} = S_{\GG'}$, i.e. that
  the $\GG'$ is purely ineffective, since $\GG$ is.

  By definition, $\GGo' = F$, and an arrow
  of $\GG'$ from $(x,\phi)$ to $(y,\psi)$ is
  an arrow $g: x \ra y$ in $\GG_{1}$ such that
  $\psi = g \phi g^{-1}$, i.e. $\psi(t) = g \phi(t) g^{-1}$ 
  for every $t \in T$. Hence the stabilizers are
  \begin{align}\label{eq:getcenter}
    S^{0}_{\GG'}(x,\phi) 
     &= \{g: x \ra x 
            \: | \: \phi(t) = g \phi(t) g^{-1} \: \: \forall t \in T \}  
        = Z(S^{0}_{\GG}(x))
  \end{align}
  since each $\phi$ is surjective.
  Define a map $c: \GGo' \times Z(T) \ra S^{0}_{\GG'}$
  by $c((x,\phi),a) = \phi(a)$. More precisely,
  we should remember the source and target of this arrow, which are both
  equal to $(x,\phi)$; hence where necessary we will denote the right hand side by
  \[
    c((x,\phi),a) = \phi(a)_{(x,\phi)}.
  \]

  By~(\ref{eq:getcenter}),
  this defines an isomorphism of local systems over $\GGo'$.
  It is $\GG'$-equivariant since, for $g:(x,\phi) \ra (y,\psi)$, 
  \[
    g c((x,\phi),a) g^{-1} = g \phi(a) g^{-1} = \psi(a) = c((y,\psi),a).
  \]
  It is also $\Aut(T)$-equivariant, where $\Aut(T)$ acts on the right
  of $Z(T)$ by $a \cdot \lambda = \lambda^{-1}(a)$.
  For, remember that $\Aut(T)$ acts on $\GGo'$ as in (\ref{autTaction}); we have
  \[
    c((x,\phi),a) \cdot \lambda = \phi(a)_{(x,\phi)} \cdot \lambda 
        = \phi(a)_{(x,\phi \circ \lambda )}
  \]
  and
  \[
    c((x,\phi) \cdot \lambda,a \cdot \lambda) = 
      c((x,\phi \circ \lambda), \lambda^{-1}(a)) = \phi(a)_{(x,\phi \circ \lambda)}. 
  \]
  It is also easy to check that $c$ is $K$-equivariant.
  Hence we have a $K$-equivariant trivialization of the stabilizer
  local system of $\GG'$, and this local system is abelian.
\end{proof}

Note: in gerbe language, this proposition accomplishes two
things. It reduces to the case of an abelian gerbe, and
it reduces to the case of a \textit{trivial band}:
\begin{Def}\label{def:trivialband}
  A gerbe $F$ on a manifold $M$, represented by a purely ineffective
  groupoid $\GG$, has \textbf{trivial abelian band} $A$, where $A$ is an abelian group, 
  if the stabilizer local system
  $S^{0}_{\GG}$ of $\GG$ is $\GG$-equivariantly 
  isomorphic to $\GGo \times A$ for some abelian group $A$.
\end{Def}
In fact one can quite generally define the notion of the band
of a gerbe, which is analogous to the coefficient sheaf for a torsor, 
but it is somewhat involved. See Giraud \cite{Giraud:Book} for the original
definition. 

Combining the two previous propositions gives the
following almost-presentation theorem.
\begin{theorem}\label{thm:presentbyPI}
  Let $\GG$ be a connected 
  orbifold groupoid of dimension $n$. Then $G$ can be presented
  as $\GG \iso \HHH/K$ where $\HHH$ is a purely ineffective orbifold groupoid
  (i.e. a manifold $P$ together with a gerbe $F$ on $P$)
  with trivial finite abelian band and an action of the compact Lie group $K$.
  The stabilizer groups of $\HHH$ are isomorphic to $Z(T)$ where
  $T$ is the ineffective stabilizer of some point of $\GG_{0}$.
\end{theorem}

As a corollary we obtain a result of Edidin-Hassett-Kresch-Vistoli 
\cite{EdidinHassettetal:Brauer}:
\begin{cor}\label{cor:trivialcenter}
  Let $\GG$ be a connected orbifold groupoid such that the ineffective
  stabilizers all have trivial center. Then $\GG$ can be presented as
  $\GG \iso P/K$ where $P$ is a smooth $K$-manifold and $K$ is a compact Lie group.
\end{cor}

We are thus led to consider the following situation. 
Let $A$ be a finite abelian group and let 
$\GG$ be a purely ineffective orbifold groupoid
such that the stabilizer local system $S^{0}_{\GG}$ of $\GG$ 
is $\GG$-equivariantly isomorphic to $\GGo \times A$.
I.e., we have a smooth manifold $M = \Gtop$ and
a gerbe $F$, with trivial abelian band.

\begin{theorem}\label{thm:presentgerbe}
  Given a smooth $m$-manifold $M$ and a gerbe $F$, 
  with trivial finite abelian band $A$. Then $(M,F)$ can be
  presented as the quotient $P/K$ of a smooth $K$-manifold $P$, where
  $K$ is a compact Lie group; in fact $K$ can be taken to be a finite
  product of unitary groups, with $A$ centrally embedded in $K$,
  and $P$ can be taken to be the total space of a principal $K/A$ bundle.
\end{theorem}

We prove a preliminary result before proving this theorem.
Let $A \iso \Z/a\Z$ be a cyclic group and let $i: A \ra U(n)$
be a central embedding of $A$ into a unitary group.
The exact sequence
\[
  \xymatrix{ 0        \ar[r] 
             & A      \ar[r]^-{i} 
             & U(n)   \ar[r]
             & U(n)/A \ar[r]
             & 0
  }
\]
gives rise to a fibration sequence 
\[
  \xymatrix{ BU(n)       \ar[r] 
             & B(U(n)/A) \ar[r]^-{\beta} 
             & K(A,2)    
  }
\]
where $\beta$ is the Bockstein map for the sequence.

For integers $n,k$, let $d_{n,k}:U(n) \ra U(n)^{k} \ra U(nk)$
be the composition of the diagonal map and the natural embedding
of $U(n)^{k}$ as block diagonal matrices in $U(n)$.
These maps form a direct system of groups and we define $U(\infty)$
to be the direct limit of this system. (Note that this is
different from $U$, which is the direct limit under direct sum
with an identity matrix.) Note that the sequence
$U(1) \ra \dots \ra U(n!) \ra U((n+1)!) \ra \dots $ is
cofinal in our direct system, so 
$U(\infty) = \displaystyle \lim_{\lra} U(n!)$.

The maps $d_{n,k}$ preserve the center, so 
given a central embedding $i:A \ra U(n)$ for some $n$, 
there will be a corresponding central embedding in any
$U(nk)$, and hence an embedding $i_{\infty}: A \ra U(\infty)$.
We have an exact sequence
\[
  0 \ra A \ra U(\infty) \ra U(\infty)/A \ra 0.
\]

The following result makes more precise a result
of Serre, cited in \cite{Grothendieck:BrauerI}.
\begin{proposition}\label{prop:surjectiveonskeleton}
  Let $A$ be a cyclic group and let
  $S \subset K(A,2)$ be an $m$-dimensional skeleton.
  Then there is an integer $n = n_{m}$ and a central embedding
  $i:A \ra U(n)$ such that the associated Bockstein 
  map $\beta: B(U(n)/A) \ra K(A,2)$ has a section over $S$. 
\end{proposition}
\begin{proof}
  We construct the section using obstruction theory.
  First let $n_{0}=1$, and let $i: A \ra U(1)$ be the standard embedding.
  We will try to section the map $\beta_{0}:B(U(n_{0})/A) \ra K(A,2)$ over
  successive skeleta of $S$, allowing ourselves to make $n_{0}$ larger,
  using the maps $d_{n,k}$ above, as follows.
  We will define a sequence 
  \[
    \xymatrix@1{ U(n_{0}) \ar[r]_{d_{n_{0},q_{0}}} & U(n_{1}) \ar[r]_{d_{n_{1},q_{1}}}
                & U(n_{2}) \ar[r]_{d_{n_{2},q_{2}}} & \dotsb }
  \]
  (where $n_{k} = q_{k-1}n_{k-1}$). Each term $U(n_{k})$ in the sequence has
  its associated map $\beta_{k}:B(U(n_{k})/A) \ra K(A,2)$, and these
  maps commute with the $d$ maps as noted above. We will inductively find a section
  of each $\beta_{k}$ over the $k$-skeleton of $S$. 

  Clearly there is a section $s_{0}$ of $\beta_{0}: B(U(n_{0})/A) \ra K(A,2)$ over
  the 0-skeleton of $K(A,2)$. Now suppose that we have 
  defined $q_{j}$ for all $j < k$ (which defines $n_{k}$)
  and we have found a section $s_{k}$ of $\beta_{k}: B(U(n_{k})/A) \ra K(A,2)$
  over the $k$-skeleton of $K(A,2)$. The obstruction to extending
  the section of $\beta_{k}$ over the $k+1$-skeleton lies in
  \[
    H^{k+1}(K(A,2);\pi_{k}(BU(n_{k}))) \iso H^{k+1}(K(A,2);\pi_{k-1}(U(n_{k}))).
  \]
  It is well-known
  that any element of $H^{k+1}(K(A,2);B)$ is torsion, for a finite abelian
  group $A$ and an arbitrary abelian group $B$. 
  So let $q_{k}$ be the least integer 
  such that every element of $H^{k+1}(K(A,2);\pi_{k-1}(U(n_{k})))$
  is $q_{k}$-torsion. So in particular the obstruction class is 
  $q_{k}$-torsion. It is easy to see 
  that under the diagonal map $d = d_{n_{k},q_{k}}$, this class dies. 
  Hence we can extend the section $d_{n_{k},q_{k}} s_{k}$ of $\beta_{k+1}$
  over the $k+1$-skeleton of $K(A,2)$. Continuing in this way, 
  we get a section of $\beta_{m}:B(U(n_{m})/A) \ra K(A,2)$
  over the $m$-skeleton, where $n_{m} = n_{0} \cdot \Prod_{k=0}^{m-1} q_{k}$.

  (Note: since $BU(n_{0})$ is 1-connected, we
  have $q_{0} = q_{1} = 1$ and so $n_{2} = n_{0}$. Higher $q_{k}$ values
  are more complicated but calculable in principle.)
\end{proof}

Now the proof of the theorem is easy.
\begin{proof}[Proof of Theorem \ref{thm:presentgerbe}]
  First we prove the theorem in the case where $A$ is cyclic.
  The gerbe $F$ is classified by an element $\alpha \in H^{2}(M;A)$.
  The class $\alpha$ defines a map $\alpha: M \ra K(A,2)$.
  Since $M$ is an $m$-manifold we can assume
  that $\alpha$ maps into an $m$-skeleton $S$ of $K(A,2)$.
  Then Proposition~\ref{prop:surjectiveonskeleton} 
  guarantees a lift of $\alpha$ through the Bockstein map
  associated to some embedding $A \ra U(n)$. This lift
  classifies a $U(n)/A$-bundle $P$ whose Bockstein is
  $\alpha$. 
  But for any such bundle we know that the class of the 
  associated $A$-gerbe is given by the Bockstein of the 
  classifying map of $P$. Hence we can present the $A$-gerbe
  using $P$.

  If $A$ is an arbitrary finite abelian group, express it
  as a product $A = \prod A_{k}$ and apply the
  above argument to each factor. We get a presentation with
  $K = \prod U(n_{k})$ and $P = P_{1} \times_{M} P_{2} \dotsb \times_{M} p_{r}$.
\end{proof}
Note that for $A$ cyclic, we can take
$n = n_{m}$ as in the proof of the last proposition. This gives,
in principle, a computable bound for how large the group $U(n)$ must
be.

In order to use this method to solve the general presentation
problem for orbifolds, one needs a positive answer to following
equivariant version of Theorem \ref{thm:presentgerbe}.

\textbf{QUESTION:}
Let $G$ be a purely ineffective orbifold groupoid
with equivariantly trivial abelian stabilizers, i.e. let
$G$ represent an abelian gerbe $F$ with finite stabilizers and
trivial band $A$ on a smooth manifold $M$. Let $K$ be a
compact Lie group and assume that $K$ acts (strictly)on the right of $G$. 
Can we find a smooth manifold $P$, a compact Lie group
$L$, and an action of $K \times L$ on $P$ such that we have
 a $K$-equivariant equivalence
\begin{equation*}
  P \rtimes L \iso G? 
\end{equation*}

\end{document}